\documentclass[reqno]{amsart}
\usepackage{geometry,amssymb,bbm}
\newcommand{\assign}{:=}
\newcommand{\dueto}[1]{\textup{\textbf{(#1) }}}
\newcommand{\mathd}{\mathrm{d}}
\newcommand{\tmop}[1]{\ensuremath{\operatorname{#1}}}
\newcommand{\tmstrong}[1]{\textbf{#1}}
\newcommand{\tmtextbf}[1]{{\bfseries{#1}}}
\newcommand{\tmtextit}[1]{{\itshape{#1}}}
\newtheorem{lemma}{Lemma}[section]
\newtheorem{remark}[lemma]{Remark}
\newtheorem{theorem}[lemma]{Theorem}

\begin{document}

\title[Global Regularity of GMHD]{Note on Solution Regularity of the Generalized Magnetohydrodynamic Equations with Partial Dissipation}
\author{Chuong V. Tran, Xinwei Yu, Zhichun Zhai}
\address{Chuong V. Tran: School of Mathematics and Statistics, 
University of St. Andrews, St Andrews KY16 9SS, United Kingdom}
\email{chuong@mcs.st-and.ac.uk}
\address{Xinwei Yu and Zhichun Zhai: Department of Mathematical and 
Statistical Sciences, University of Alberta, Edmonton, AB, T6G 2G1, Canada}
\email{xinweiyu@math.ualberta.ca, zhichun1@ualberta.ca}



\begin{abstract}
  In this brief note we study the $n$-dimensional magnetohydrodynamic 
  equations with hyper-viscosity and zero resistivity. We prove global
  regularity of solutions when the hyper-viscosity is sufficiently strong.
\end{abstract}
\maketitle

\section{Introduction}

Consider the $n$-dimensional generalized magnetohydrodynamic ($n$D GMHD)
equations 
\begin{eqnarray}
  u_t + u \cdot \nabla u & = & - \nabla p + b \cdot \nabla b - \nu
  \mathcal{L}_1^2 u,  \label{eq:GMHD-u}\\
  b_t + u \cdot \nabla b & = & b \cdot \nabla u - \kappa \mathcal{L}_2^2 b, 
  \label{eq:GMHD-b}\\
  \nabla \cdot u = \nabla \cdot b & = & 0,  \label{eq:GMHD-div-free}
\end{eqnarray}
where the Laplacians $\triangle$ in the dissipation terms of the momentum
and induction equations have been replaced by general negative-definite 
operators $-\mathcal{L}_1^2$ and $-\mathcal{L}_2^2$, respectively. Various 
forms of these operators have been used in studies concerning the persistence of regularity for classical solutions. In particular, Wu \cite{Wu2003} considered 
$\mathcal{L}_1 = \Lambda^{\alpha}$ and $\mathcal{L}_2 = \Lambda^{\beta}$, 
where $\Lambda \assign \left( - \triangle \right)^{1 / 2}$, and proved 
global regularity, that is classical solutions exist for all time, when both $\alpha \geqslant \frac{1}{2} + \frac{n}{4}$ 
and $\beta \geqslant \frac{1}{2} + \frac{n}{4}$ hold concurrently. This 
result has been improved by several authors 
\cite{Wu2009,Wu2008,Wu2011,Zhou2007} (also see \cite{Cao2011} for the 
case of degenerate $\mathcal{L}_i$'s). To date, the best global regularity 
result for (\ref{eq:GMHD-u}--\ref{eq:GMHD-div-free}) is the following theorem.\

{\noindent}\tmtextbf{Theorem. }\tmtextit{{\dueto{Wu 2011
{\cite{Wu2011}}}}Consider the GMHD system
(\ref{eq:GMHD-u}--\ref{eq:GMHD-div-free}) with $\mathcal{L}_1, \mathcal{L}_2$
defined through Fourier transform as
\begin{equation}
  \widehat{\mathcal{L}_1 u} \left( \xi \right) = m_1 \left( \xi \right) 
  \widehat{u} \left( \xi \right), \hspace{2em} \widehat{\mathcal{L}_2 b} \left(
  \xi \right) = m_2 \left( \xi \right)  \widehat{b} \left( \xi \right)
\end{equation}
with
\begin{equation}
  m_1 \left( \xi \right) \geqslant \frac{\left| \xi \right|^{\alpha}}{g_1
  \left( \left| \xi \right| \right)}, \hspace{2em} m_2 \left( \xi \right)
  \geqslant \frac{\left| \xi \right|^{\beta}}{g_2 \left( \left| \xi \right|
  \right)}
\end{equation}
where $g_1 \geqslant 1$ and $g_2 \geqslant 1$ are nondecreasing. Assume the initial data belong to $H^s$ with $s>1+\frac{n}{2}$. Then 
the system has a unique global classical solution if the following conditions are satisfied:
\begin{equation}
  \alpha \geqslant \frac{1}{2} + \frac{n}{4}, \hspace{2em} \beta > 0,
  \hspace{2em} \alpha + \beta \geqslant 1 + \frac{n}{2}, \hspace{2em}
  \int_1^{\infty} \frac{\mathd s}{s \left( g_1 \left( s \right)^2 + g_2 \left(
  s \right)^2 \right)^2} = + \infty . \label{eq:Wu2011-cond}
\end{equation}}{\hspace*{\fill}}{\medskip}

When $n = 2$, conditions much weaker than (\ref{eq:Wu2011-cond}) are
sufficient \cite{Tran2012}. For example, in the absence of viscosity 
(i.e. $\nu=0$), global regularity can be secured provided $\beta>2$ 
(and $g_2=1$). For $n \geqslant 3$, such a complete removal of  
$\mathcal{L}_1$ is inconceivable. In fact, a drastic weakening of 
$\mathcal{L}_1$ can hardly be expected. The reason is that the 
equations (\ref{eq:GMHD-u}--\ref{eq:GMHD-div-free}) contain the
generalized Navier-Stokes system
\begin{equation}
  u_t + u \cdot \nabla u = - \nabla p - \frac{\Lambda^{2 \alpha}}{g_1 \left(
  \Lambda \right)^2} u, \hspace{2em} \nabla \cdot u = 0
\end{equation}
as a special case (obtained by setting $b = 0$), for which the problem 
of global regularity is still open unless \cite{Tao2009}
\begin{equation}
  \alpha \geqslant \frac{1}{2} + \frac{n}{4}, \hspace{2em}
  \int_1^{\infty} \frac{\mathd s}{s g_1(s)^4} = + \infty . 
\label{eq:Tao2009-cond}
\end{equation}
(also see \cite{Zhang2010} for the anisotropic case) Hence, before (\ref{eq:Tao2009-cond}), which consists of the first 
and final conditions (in the absence of $g_2$) in (\ref{eq:Wu2011-cond}), 
can be weakened, an improvement concerning the conditions on $\alpha$ 
and $g_1$ in Wu's theorem is highly infeasible.

As indicated by the discussion in the preceding paragraph, the condition 
on $\alpha$ in (\ref{eq:Wu2011-cond}) is ``genuine,'' and its weakening 
would be a formidable task. On the other hand, the condition $\beta>0$ 
appears ``technical'' and could be removed. The intuitive reason is 
that for a sufficiently strong $\mathcal{L}_1$, bounds can be derived 
for sufficiently high order derivatives of $u$. Since the induction
equation is linear in $b$, this result in turn can be used to prove 
boundedness for sufficiently high order derivatives 
of $b$, even in the absence of magnetic diffusion, thereby ensuring 
regularity. The question is whether the removal of $\beta>0$ can be 
done without a cost. It turns out that the answer to this question is 
positive. In fact, we show in this article that 
(\ref{eq:Wu2011-cond}) can be readily extended to the case $\beta=0$ 
(more precisely to $\kappa=0$). This is accomplished through an 
application of Lei and Zhou's ``weakly nonlinear'' energy estimate 
approach \cite{Lei2009}, which enables us to derive ``almost a priori'' 
bounds for the $H^1$ norms of $u$ and $b$. These results are sufficient 
for obtaining uniform bounds for higher Sobolev norms, hence implying
global regularity. To the best of our knowledge, Lei and Zhou first 
applied this approach to mathematical fluid mechanics in {\cite{Lei2009}}.

Now we state our main result.

\begin{theorem}
  \label{thm:main}Consider the following GMHD system
  \begin{eqnarray}
    u_t + u \cdot \nabla u & = & - \nabla p + b \cdot \nabla b - \nu
    \mathcal{L}^2 u,  \label{eq:GMHD-u-thm}\\
    b_t + u \cdot \nabla b & = & b \cdot \nabla u,  \label{eq:GMHD-b-thm}\\
    \nabla \cdot u = \nabla \cdot b & = & 0,  \label{eq:GMHD-div-free-thm}
  \end{eqnarray}
  with
  \begin{equation}
    \mathcal{L} \assign \frac{\Lambda^{\alpha}}{g \left( \Lambda \right)}
    \hspace{2em} \text{\tmop{defined} \tmop{as}} ~ \widehat{\mathcal{L}u} 
    \left( \xi \right) \assign \frac{\left| \xi \right|^{\alpha}}{g 
    \left( \left| \xi \right| \right)}  \widehat{u} \left( \xi \right),
  \end{equation}
  for some function $g \left( s \right) \geqslant 1$ defined on $s \geqslant
  0$. Let the initial data $u_0, b_0 \in H^k$ for some $k > 1 + \frac{n}{2}$.
  Then the system has a unique global classical solution if the following conditions are
  satisfied:
  \begin{equation}
    \alpha \geqslant 1 + \frac{n}{2}, \hspace{2em} g \left( s \right)^2
    \leqslant C \log \left( e + s \right) ~ \text{\tmop{for} \tmop{some}
    \tmop{absolute} \tmop{constant}} ~ C. \label{eq:cond-main}
  \end{equation}
\end{theorem}

The following remarks are in order.
\begin{itemize}
  \item It is clear that (\ref{eq:cond-main}) extends (\ref{eq:Wu2011-cond})
  to the case $\beta = 0$.

  \item $g(s)$ does not need to be nondecreasing.
  
  \item In some sense, the condition $g \left( s \right)^2 \leqslant C 
  \log \left( e + s \right)$ is weaker than $\int_1^{\infty} \frac{\mathd
  s}{sg \left( s \right)^4} = + \infty$. For example, given the typical 
  case $g \left( s \right) \sim \left[ \log \left( e + s \right)
  \right]^{\gamma}$, the former requires $\gamma \leqslant 1 / 2$ while 
   the latter requires $\gamma \leqslant 1 / 4$. 
\end{itemize}

The remaining of this article is devoted to the proof of Theorem
\ref{thm:main}. In what follows, we set $\nu = 1$ to simplify the 
presentation. The adaptation of the proof for other values of $\nu$ 
is straightforward.

\section{Proof of Theorem \ref{thm:main}}

We present detailed proof for the case $\alpha = 1 + \frac{n}{2}$, that
is $\mathcal{L} \assign \frac{\Lambda^{1 + n / 2}}{g \left( \Lambda \right)}$,
or more explicitly
\begin{equation}
  \widehat{\mathcal{L}u} \left( \xi \right) = \frac{\left| \xi \right|^{1 + n
  / 2}}{g \left( \left|\xi\right| \right)}  \widehat{u} \left( \xi \right) .
\end{equation}
The case $\alpha > 1 + \frac{n}{2}$ is much easier to handle and 
briefly discussed at the end of this section.

Multiplying (\ref{eq:GMHD-u-thm}) and (\ref{eq:GMHD-b-thm}) by $u$ and 
$b$, respectively, and integrating the resulting equations over space, 
we obtain the standard energy equality
\begin{equation}
  \frac{\mathd}{\mathd t} \frac{\left\| u \right\|_{L^2}^2 + \left\| b
  \right\|_{L^2}^2}{2} + \left\| \mathcal{L}u \right\|_{L^2}^2 = 0.
\label{eq:energy}
\end{equation}
Integrating (\ref{eq:energy}) up to some fixed (but arbitrary) time $T$, we deduce that
\begin{equation}
  u, b \in L^{\infty} \left( 0, T ; L^2 \right), \hspace{2em} \mathcal{L}u \in
  L^2 \left( 0, T ; L^2 \right) . \label{eq:energy-est}
\end{equation}
Note that for any $0 \leqslant \lambda < 1 + n / 2$ and any $m\geqslant 0$ , there is a constant $C$
depending only on $\lambda$, $m$, and $g$ such that
\begin{equation}
  \left\| u \right\|_{H^{m + \lambda}} \leqslant C \left( \left\| u
  \right\|_{L^2} + \left\| \mathcal{L} \Lambda^m u \right\|_{L^2} \right).
  \label{eq:L2lambda}
\end{equation}
This, together with (\ref{eq:energy-est}), implies that 
$u \in L^2 \left( 0, T ; H^{\lambda} \right)$ for any 
$0 \leqslant \lambda < 1 + n / 2$.

In the following, we will show that for any $T > 0$, 
$\left\| u \right\|_{H^k}$ and $\left\| b \right\|_{H^k}$
are uniformly bounded over $\left( 0, T \right)$, or more
precisely over $\left( T_0, T \right)$ for some $T_0$ close enough to $T$. As
local well-posedness for (\ref{eq:GMHD-u-thm}--\ref{eq:GMHD-div-free-thm}) can
be proved by standard methods, such uniform bounds secure global regularity. We first show that under the assumption of
Theorem \ref{thm:main}, once (\ref{eq:energy-est}) holds, the $H^1$ norms of
$u, b$ have to be much smaller than their $H^k$ norms. This makes the
trilinear terms in the standard energy method much weaker than its scaling
suggests, thereby enabling us to derive $H^k$ a priori bounds.

\subsection{$H^1$ Estimates}

The key to our derivation of estimates in $H^1$ is the following lemma, 
whose proof is given in the appendix.

\begin{lemma}
  \label{lem:log-ineq}Let $g : \mathbbm{R}^+ \mapsto \left[ 1, + \infty
  \right)$ be such that $g \left( s \right)^2 \leqslant C_0 \log \left( e + s
  \right)$ for some absolute constant $C_0$ and for all $s \geqslant 0$, then there is a constant $C=C(k,n)$ such that
  \begin{equation}
    \left\| \nabla u \right\|_{L^{\infty}} \leqslant C \left[ \left\| u \right\|_{L^2}
    + \left\| \frac{\Lambda^{1 + \frac{n}{2}}}{g \left( \Lambda \right)} u
    \right\|_{L^2} \log \left( e + \left\| u \right\|_{H^k} \right) \right]
  \end{equation}
  for any $k > 1 + \frac{n}{2}$. 
\end{lemma}

\begin{remark}
  This can be seen as a variant of the classical Brezis-Wainger inequality
  (see e.g. {\cite{Brezis1980}}, {\cite{Engler1989}}) where $g \left( \Lambda
  \right) = 1$ and the $\log$ factor is $\left( \log \left( e + \left\| u
  \right\|_{H^k} \right) \right)^{1 / 2}$. It can also be seen as a limiting case of the Sobolev inequalities (see e.g. {\cite{Kozono2000}}).
\end{remark}

Let $\partial_i$ denote a partial derivative. Differentiating 
(\ref{eq:GMHD-u-thm},\ref{eq:GMHD-b-thm}) yields 
\begin{eqnarray}
  \label{eq:partial1}
  \left( \partial_i u \right)_t + u \cdot \nabla \partial_i u & = & -
  \partial_i u \cdot \nabla u - \nabla \partial_i p + \partial_i b \cdot
  \nabla b + b \cdot \nabla \partial_i b -\mathcal{L}^2 \partial_i u, \\
  \label{eq:partial2}
  \left( \partial_i b \right)_t + u \cdot \nabla \partial_i b & = & -
  \partial_i u \cdot \nabla b + \partial_i b \cdot \nabla u + b \cdot \nabla
  \partial_i u. 
\end{eqnarray}
Multiplying (\ref{eq:partial1}) and (\ref{eq:partial2}) by $\partial_i u$ 
and $\partial_i b$, respectively, integrating the resulting equations in 
space and summing up over $i$ (noting $\nabla \cdot u = \nabla \cdot b = 0$) 
we obtain
\begin{equation}
  \frac{\mathd}{\mathd t} \int_{\mathbbm{R}^n} \left( \frac{\left| \nabla u
  \right|^2 + \left| \nabla b \right|^2}{2} \right) \mathd x +   \int_{\mathbbm{R}^n} \left| \mathcal{L} \nabla u \right|^2 \mathd x
  \leqslant C \left\| \nabla u \right\|_{L^{\infty}}  \int_{\mathbbm{R}^n}
  \left( \left| \nabla u \right|^2 + \left| \nabla b \right|^2 \right) \mathd
  x.
\end{equation}
This implies
\begin{equation}
  \left( \left\| \nabla u \right\|_{L^2}^2 + \left\| \nabla b \right\|_{L^2}^2
  \right) \left( t \right) \leqslant \left( \left\| \nabla u \right\|_{L^2}^2
  + \left\| \nabla b \right\|_{L^2}^2 \right) \left( T_0 \right) \exp \left[
  C\int_{T_0}^t \left\| \nabla u \right\|_{L^{\infty}} \left( \tau \right)
  \mathd \tau \right] .
\end{equation}
Applying Lemma \ref{lem:log-ineq} we have for any $T_0 < t$,
\begin{eqnarray}
  \left( \left\| \nabla u \right\|_{L^2}^2 + \left\| \nabla b \right\|_{L^2}^2
  \right) \left( t \right) & \leqslant & C \left( T_0 \right) \exp \left[
  \int_{T_0}^t \left( \left\| u \right\|_{L^2} + C \left\| \mathcal{L}u
  \right\|_{L^2} \log \left( e + \left\| u \right\|_{H^k} \right) \right)
  \mathd s \right] \nonumber\\
  & \leqslant & C \left( T_0 \right) \exp \left[ C \left( \int_{T_0}^t
  \left\| \mathcal{L}u \right\|_{L^2} \mathd \tau \right) \log \left( M \left(
  t \right) \right) \right] \nonumber\\
  & \leqslant & C \left( T_0 \right) M \left( t \right)^{C \left(
  \int_{T_0}^t \left\| \mathcal{L}u \right\|_{L^2} \mathd \tau \right)} ,
\end{eqnarray}
where
\begin{equation}
  M \left( t \right) \assign \sup_{\tau \in \left( T_0, t \right)} \left[ e +
  \left\| u \right\|_{H^k} + \left\| b \right\|_{H^k} \right] .
  \label{eq:def-M}
\end{equation}
Note that we have used $\left\| u \right\|_{L^2} \leqslant
\left\| u_0 \right\|_{L^2} + \left\| b_0 \right\|_{L^2}$. Also 
note that the value of $C\left( T_0 \right)$ changes from line to line.

As $k > 1 + \frac{n}{2}$, there exists $\lambda$ satisfying
\begin{equation}
  \frac{n}{2}  \frac{k}{k - 1} < \lambda < 1 + \frac{n}{2} .
  \label{eq:choice-lambda}
\end{equation}
Now using $\left\| \mathcal{L}u \right\|_{L^2} \in L^2 \left( 0, T \right)$ we
see that there exists $T_0 < T$ such that for all $t \in \left( T_0, T
\right)$,
\begin{equation}
  C \int_{T_0}^t \left\| \mathcal{L}u \right\|_{L^2} \mathd \tau < 2 \delta
  \assign \min \left( \frac{\left( k + \lambda \right)  \left( k - 1 \right) -
  k \left( k - 1 + \frac{n}{2} \right)}{k \left( k - 1 - \frac{n}{2} \right)},
  \frac{\lambda - \frac{n}{2}}{k + \lambda} \right) . \label{eq:choice-delta}
\end{equation}
Thanks to (\ref{eq:choice-lambda}), the right-hand side of 
(\ref{eq:choice-delta}) is positive since both numbers in the 
brackets are positive. This allows us to fix $T_0$. In what folows,
$T_0$ is thus fixed.

\subsection{$H^k$ Estimates}

Let $\partial^k$ denote any $k$th order partial derivative. By applying
$\partial^k$ to each of (\ref{eq:GMHD-u-thm}) and (\ref{eq:GMHD-b-thm}),
multiplying the resulting equations by $\partial^k u$ and $\partial^k b$, 
respectively, and integrating we obtain
\begin{eqnarray}
  \frac{\mathd}{\mathd t} \left( \frac{\left\| \partial^k u \right\|_{L^2}^2 +
  \left\| \partial^k b \right\|_{L^2}^2}{2} \right) + \left\| \mathcal{L}
  \partial^k u \right\|_{L^2} & = & - \int_{\mathbbm{R}^n} \partial^k \left( u
  \cdot \nabla u \right) \partial^k u \mathd x + \int_{\mathbbm{R}^n}
  \partial^k \left( b \cdot \nabla b \right) \partial^k u \mathd x \nonumber\\
  &  & - \int_{\mathbbm{R}^n} \partial^k \left( u \cdot \nabla b \right)
  \partial^k b \mathd x + \int_{\mathbbm{R}^n} \partial^k \left( b \cdot
  \nabla u \right) \partial^k b \mathd x. 
\end{eqnarray}
Now summing over all $k$th partial derivatives, and taking advantage of
$\nabla \cdot u = \nabla \cdot b = 0$, we reach
\begin{equation}
  \frac{\mathd}{\mathd t} \left(\frac{ \left\| \nabla^k u \right\|_{L^2}^2 + \left\|
  \nabla^k b \right\|_{L^2}^2}{2} \right) +  \left\| \mathcal{L} \nabla^k u
  \right\|_{L^2}^2 = I_1 + I_2 + I_3 ,
\end{equation}
where
\begin{eqnarray}
  I_1 & = & - \sum \int_{\mathbbm{R}^n} \left[ \partial^k \left( u \cdot
  \nabla u \right) - u \cdot \nabla \partial^k u \right] \partial^k u \mathd
  x, \\
  I_2 & = & \sum \int_{\mathbbm{R}^n} \left[ \partial^k \left( b \cdot \nabla
  b \right) - b \cdot \nabla \partial^k b \right] \partial^k u \mathd x +
  \int_{\mathbbm{R}^n} \left[ \partial^k \left( b \cdot \nabla u \right) - b
  \cdot \nabla \partial^k u \right] \partial^k b \mathd x, \\
  I_3 & = & - \sum \int_{\mathbbm{R}^n} \left[ \partial^k \left( u \cdot
  \nabla b \right) - u \cdot \nabla \partial^k b \right] \partial^k b \mathd
  x. 
\end{eqnarray}
From this we see that
\begin{eqnarray}
  \frac{\mathd}{\mathd t} \left( \frac{\left\| \nabla^k u \right\|_{L^2}^2 + \left\|
  \nabla^k b \right\|_{L^2}^2}{2} \right) +  \left\| \mathcal{L} \nabla^k u
  \right\|_{L^2}^2 & \leqslant & \sum \left| \int_{\mathbbm{R}^n} \partial^l u
  \partial^m u \partial^k u \mathd x \right| + \sum \left|
  \int_{\mathbbm{R}^n} \partial^l b \partial^m b \partial^k u \mathd x \right|
  \nonumber\\
  &  & + \sum \left| \int_{\mathbbm{R}^n} \partial^l b \partial^m u
  \partial^k b \mathd x \right| . 
\end{eqnarray}
The summation is over all possible combinations of partial derivatives
satisfying $l + m = k + 1, l, m \geqslant 1$.
\begin{itemize}
  \item Estimating $\sum \left| \int_{\mathbbm{R}^n} \partial^l u \partial^m u
  \partial^k u \mathd x \right| + \sum \left| \int_{\mathbbm{R}^n} \partial^l
  b \partial^m b \partial^k u \mathd x \right|$.
  
  These terms can be estimated similarly. So we only present detailed 
  calculations for 
  $\left| \int_{\mathbbm{R}^n} \partial^l b \partial^m b \partial^k u 
  \mathd x \right|$. 
  
  First applying H\"oder's inequality to the integral yields
  \begin{equation}
    \left| \int_{\mathbbm{R}^n} \partial^l b \partial^m b \partial^k u \mathd
    x \right| \leqslant \left\| \partial^l b \right\|_{L^2}  \left\| \partial^m b
    \right\|_{L^2}  \left\| \partial^k u \right\|_{\infty} .
  \end{equation}
  Thanks to (\ref{eq:choice-lambda}) the following Gagliardo-Nirenberg
  inequality holds:
  \begin{equation}
    \left\| \partial^k u \right\|_{L^{\infty}} \leqslant C \left\| u
    \right\|_{L^2}^a \left\| \Lambda^{k + \lambda} u \right\|_{L^2}^{1 - a}
  \end{equation}
  with
  \begin{equation}
    a = \frac{\lambda - \frac{n}{2}}{k + \lambda} \Longrightarrow 1 - a =
    \frac{k + \frac{n}{2}}{k + \lambda} . \label{eq:a}
  \end{equation}
  Furthermore, as $l, m \geqslant 1$, we have
  \begin{equation}
    \left\| \partial^l b \right\|_{L^2} \leqslant C \left\| \nabla\ b
    \right\|_{L^2}^{\xi}  \left\| \nabla^k b \right\|_{L^2}^{1 - \xi} ;
    \hspace{2em} \left\| \partial^m b \right\|_{L^2} \leqslant C \left\| \nabla
    b \right\|_{L^2}^{\eta}  \left\| \nabla^k b \right\|_{L^2}^{1 - \eta}
  \end{equation}
  with
  \begin{equation}
    \xi = \frac{k - l}{k - 1}, \hspace{1em} \eta = \frac{k - m}{k - 1} .
  \end{equation}
  Thus we reach
  \begin{equation}
    \left| \int_{\mathbbm{R}^n} \partial^l b \partial^m b \partial^k u \mathd
    x \right| \leqslant C \left\| \nabla b \right\|_{L^2}  \left\|
    \nabla^k b \right\|_{L^2}  \left\| u \right\|_{L^2}^a  \left\| \Lambda^{k +
    \lambda} u \right\|_{L^2}^{1 - a}.
  \end{equation}
  As $a > 0$ we have $1 + 1 - a < 2$ and therefore can apply Young's
  inequality to get
  \begin{equation}
    \left| \int_{\mathbbm{R}^n} \partial^l b \partial^m b \partial^k u \mathd
    x \right| \leqslant C \left\| \nabla b \right\|_{L^2}^{\frac{2}{1 + a}} 
    \left\| \nabla^k b \right\|_{L^2}^{\frac{2}{1 + a}}  \left\| u
    \right\|_{L^2}^{\frac{2 a}{1 + a}} + \varepsilon \left\| \Lambda^{k +
    \lambda} u \right\|_{L^2}^2 .
  \end{equation}
  Now using $\left\| u \right\|_{L^2}\leqslant \left\| u_0 \right\|_{L^2} + \left\| b_0 \right\|_{L^2}$ and
  (\ref{eq:L2lambda}) we conclude that
  \begin{equation}
    \left| \int_{\mathbbm{R}^n} \partial^l b \partial^m b \partial^k u \mathd
    x \right| \leqslant C \left\| \nabla b \right\|_{L^2}^{\frac{2}{1 + a}} 
    \left\| \nabla^k b \right\|_{L^2}^{\frac{2}{1 + a}} + \varepsilon \left[
    \left\| \mathcal{L}\Lambda^k u \right\|_{L^2}^2 + 1 \right]
  \end{equation}
  for $\varepsilon$ as small as necessary.
  
  As the other term can be estimated similarly, we obtain, after taking
  an appropriate value of $\varepsilon$,
  \begin{eqnarray}
    \sum \left[ \left| \int_{\mathbbm{R}^n} \partial^l u \partial^m u \partial^k u
    \mathd x \right| + \left| \int_{\mathbbm{R}^n} \partial^l b \partial^m b
    \partial^k u \mathd x \right|\right] & \leqslant & C \left[ \left\| \nabla u
    \right\|_{L^2}^{\frac{2}{1 + a}}  \left\| \nabla^k u
    \right\|_{L^2}^{\frac{2}{1 + a}} + \left\| \nabla b
    \right\|_{L^2}^{\frac{2}{1 + a}}  \left\| \nabla^k b
    \right\|_{L^2}^{\frac{2}{1 + a}} \right] \nonumber\\
    &  & + \frac{1}{4}  \left[ \left\| \mathcal{L} \Lambda^k u \right\|_{L^2}^2
    + 1 \right].  \label{eq:est-1st-2nd}
  \end{eqnarray}
  \begin{remark}
    Note that the $H^1$ estimates play crucial roles here. Without them 
    we would have to use \
    \begin{equation}
      \left\| \partial^l b \right\|_{L^2} \leqslant \left\| b
      \right\|_{L^2}^{\xi}  \left\| \nabla^k b \right\|_{L^2}^{1 - \xi},
      \hspace{2em} \left\| \partial^m b \right\|_{L^2} \leqslant C \left\| b
      \right\|_{L^2}^{\eta}  \left\| \nabla^k b \right\|_{L^2}^{1 - \eta}
    \end{equation}
    with
    \begin{equation}
      1 - \xi = \frac{l}{k}, \hspace{2em} 1 - \eta = \frac{m}{k} 
    \end{equation}
    and end up with
    \begin{equation}
      \left| \int_{\mathbbm{R}^n} \partial^l b \partial^m b \partial^k u
      \mathd x \right| \leqslant C \left\| \nabla^k b \right\|^{\frac{l +
      m}{k}}_{L^2}  \left\| \Lambda^{k + \lambda} u \right\|_{L^2}^{1 - a} .
    \end{equation}
    Now applying Young's inequality would yield the term
    $\|\mathcal{L}\Lambda^k u\|^\gamma$, where $\gamma>2$, because
    \begin{equation}
      \frac{l + m}{k} + 1 - a = \frac{k + 1}{k} + 1 - a = 2 + \frac{1}{k} -
      \frac{\lambda - \frac{n}{2}}{k + \lambda} = 2 + \frac{k + \lambda -
      \left( \lambda - \frac{n}{2} \right) k}{k \left( k + \lambda \right)} >
      2
    \end{equation}
    for all $\lambda < 1 + \frac{n}{2}$. Apparently, such a term is beyond 
    the control of the available dissipation term.
  \end{remark}
  
  \item Estimating $\sum \left| \int_{\mathbbm{R}^n} \partial^l b \partial^m u
  \partial^k b \mathd x \right|$.
  
  We first apply H\"older's inequality
  \begin{equation}
    \left| \int_{\mathbbm{R}^n} \partial^l b \partial^m u \partial^k b \mathd
    x \right| \leqslant \left\| \partial^l b \partial^m u \right\|_{L^2} 
    \left\| \partial^k b \right\|_{L^2} .
  \end{equation}
  Now the standard calculus inequality (see e.g. {\cite{Majda2002}}) gives
  (recall that $l, m \geqslant 1$):
  \begin{equation}
    \left| \int_{\mathbbm{R}^n} \partial^l b \partial^m u \partial^k b \mathd
    x \right| \leqslant C \left[ \left\| \nabla u \right\|_{L^{\infty}} 
    \left\| \nabla^k b \right\|_{L^2}^2 + \left\| \nabla b
    \right\|_{L^{\infty}}  \left\| \nabla^k u \right\|_{L^2}  \left\| \nabla^k
    b \right\|_{L^2} \right] .
  \end{equation}
  For the first term on the right-hand side, applying Lemma 
  \ref{lem:log-ineq} yields
  \begin{equation}
    \left\| \nabla u \right\|_{L^{\infty}}  \left\| \nabla^k b \right\|_{L^2}
    \leqslant C \left[ 1 + \left\| \mathcal{L}u \right\|_{L^2} \log \left( e +
    \left\| u \right\|_{H^k} + \left\| b \right\|_{H^k} \right) \right] 
    \left\| \nabla^k b \right\|_{L^2} .
  \end{equation}
  For the second term, we resort to the following Gagliardo-Nirenberg 
  inequalities. First, we have
  \begin{equation}
    \left\| \nabla b \right\|_{\infty} \leqslant C \left\| \nabla b
    \right\|_{L^2}^{\xi}  \left\| \nabla^k b \right\|_{L^2}^{1 - \xi},
  \end{equation}
  where
  \begin{equation}
    \xi = \frac{k - 1 - \frac{n}{2}}{k - 1} \Longrightarrow 1 - \xi = \frac{n
    / 2}{k - 1} .
  \end{equation}
  Second, 
  \begin{equation}
    \left\| \nabla^k u \right\|_{L^2} \leqslant C \left\| \Lambda^{\lambda} u
    \right\|_{L^2}^{\eta}  \left\| \Lambda^{k + \lambda} u \right\|_{L^2}^{1 -
    \eta} ,
  \end{equation}
  where
  \begin{equation}
    \eta = \frac{\lambda}{k} \Longrightarrow 1 - \eta = \frac{k - \lambda}{k} .
  \end{equation}
  Note that $\lambda < 1 + \frac{n}{2} < k$. It follows that
  \begin{equation}
    \left\| \nabla b \right\|_{L^{\infty}}  \left\| \nabla^k u \right\|_{L^2} 
    \left\| \nabla^k b \right\|_{L^2} \leqslant C \left\| \nabla b
    \right\|_{L^2}^{\xi}  \left\| \nabla^k b \right\|_{L^2}^{2 - \xi}  \left\|
    \Lambda^{\lambda} u \right\|_{L^2}^{\eta}  \left\| \Lambda^{k + \lambda} u
    \right\|_{L^2}^{1 - \eta} .
  \end{equation}
 Obviously $\xi + \eta \leqslant 2$. Furthermore, thanks to  
  (\ref{eq:choice-lambda}), we have
  \begin{equation}
    \xi + \eta = \frac{k - 1 - \frac{n}{2}}{k - 1} + \frac{\lambda}{k} = 1 +
    \frac{\lambda}{k} - \frac{\frac{n}{2}}{k - 1} = 1 + \frac{ \left( k - 1
    \right) \lambda - \frac{n}{2} k}{k \left( k - 1 \right)} > 1 .
  \end{equation}
 Therefore
  \begin{equation}
    2 - \xi + 1 - \eta < 2, \hspace{2em} 2 - \xi \geqslant \eta .
  \end{equation}
 This enables us to apply Young's inequality to obtain
  \begin{equation}
    \left\| \nabla b \right\|_{L^{\infty}}  \left\| \nabla^k u \right\|_{L^2} 
    \left\| \nabla^k b \right\|_{L^2} \leqslant C \left\| \nabla b \right\|^A 
    \left\| \nabla^k b \right\|^B  \left\| \Lambda^{\lambda} u \right\|^C +
    \varepsilon \left\| \Lambda^{k + \lambda} u \right\|^2
  \end{equation}
  with
  \begin{equation}
    A = \frac{2 k \left( k - 1 - \frac{n}{2} \right)}{\left( k + \lambda
    \right)  \left( k - 1 \right)}, \hspace{1em} B = \frac{2 k \left( k - 1 +
    \frac{n}{2} \right)}{\left( k + \lambda \right)  \left( k - 1 \right)} <
    2, \hspace{1em} C = \frac{2 \lambda}{k + \lambda} < 2 .\label{eq:ABC}
  \end{equation}
  Now by (\ref{eq:L2lambda}) and $\left\| u \right\|_{L^2}
  \leqslant \left\| u_0 \right\|_{L^2} + \left\| b_0 \right\|_{L^2}$, we have
  $\left\| \Lambda^{k + \lambda} u \right\|_{L^2} \leqslant C \left( \left\| u
  \right\|_{L^2} + \left\| \mathcal{L} \Lambda^k u \right\|_{L^2} \right)
  \leqslant C \left( 1 + \left\| \mathcal{L} \Lambda^k u \right\|_{L^2}
  \right)$. Similarly, $\left\| \Lambda^{\lambda} u \right\|_{L^2} \leqslant
  C \left( 1 + \left\| \mathcal{L}u \right\|_{L^2} \right)$. So finally we
  reach
  \begin{eqnarray}
    & \left| \int_{\mathbbm{R}^n} \partial^l b \partial^m u \partial^k b
    \mathd x \right| \leqslant & C \left[ 1 + \left\| \mathcal{L}u
    \right\|_{L^2} \log \left( e + \left\| u \right\|_{H^k} + \left\| b
    \right\|_{H^k} \right) \right]  \left\| \nabla^k b \right\|_{L^2}
    \nonumber\\
    &  & + C \left\| \nabla b \right\|^A  \left\| \nabla^k b \right\|^B 
    \left( 1 + \left\| \mathcal{L}u \right\|_{L^2} \right) + \frac{1}{4} 
    \left( 1 + \left\| \mathcal{L} \Lambda^k u \right\|_{L^2}^2 \right) . 
  \end{eqnarray}
\end{itemize}
In summary, we have obtained
\begin{eqnarray}
  \frac{\mathd}{\mathd t} \left( \frac{\left\| \nabla^k u \right\|_{L^2}^2 + \left\|
  \nabla^k b \right\|_{L^2}^2}{2} \right) + \left\| \mathcal{L} \nabla^k u
  \right\|_{L^2}^2 & \leqslant & C \left[ \left\| \nabla u
  \right\|_{L^2}^{\frac{2}{1 + a}}  \left\| \nabla^k u
  \right\|_{L^2}^{\frac{2}{1 + a}} + \left\| \nabla b
  \right\|_{L^2}^{\frac{2}{1 + a}}  \left\| \nabla^k b
  \right\|_{L^2}^{\frac{2}{1 + a}} \right] \nonumber\\
  && + \frac{1}{4} \left(1+\left\| \mathcal{L}\Lambda^k u \right\|_{L^2}^2 \right) \nonumber\\
  &  & + C \left[ 1 + \left\| \mathcal{L}u \right\|_{L^2} \log \left( e +
  \left\| u \right\|_{H^k} + \left\| b \right\|_{H^k} \right) \right]  \left\|
  \nabla^k b \right\|_{L^2}^2 \nonumber\\
  &  & + C \left\| \nabla b \right\|^A  \left\| \nabla^k b \right\|^B  \left(
  1 + \left\| \mathcal{L}u \right\|_{L^2} \right) \nonumber\\
  &  & + \frac{1}{4}  \left( 1 + \left\| \mathcal{L} \Lambda^k u
  \right\|_{L^2}^2 \right) . 
\end{eqnarray}
Here $A, B$ and $a$ are defined in (\ref{eq:a}) and (\ref{eq:ABC}). Recalling
the definition of $M \left( t \right)$ in (\ref{eq:def-M}), we have
\begin{eqnarray}
  \frac{\mathd}{\mathd t} \left( \left\| \nabla^k u \right\|_{L^2}^2 + \left\|
  \nabla^k b \right\|_{L^2}^2 \right) & \leqslant & C \left( \left\| \nabla u
  \right\|_{L^2}^{\frac{2}{1 + a}} + \left\| \nabla b
  \right\|_{L^2}^{\frac{2}{1 + a}}  \right) M \left( t \right)^{\frac{2}{1 +
  a}} \nonumber\\
  &  & + C \left[ 1 + \left\| \mathcal{L}u \right\|_{L^2} \log \left( M
  \left( t \right) \right) \right] M \left( t \right)^2 \nonumber\\
  &  & + C \left\| \nabla b \right\|^A M \left( t \right)^B  \left( 1 +
  \left\| \mathcal{L}u \right\|_{L^2} \right) . 
\end{eqnarray}
Here we have used the fact that by definition $M \left( t \right)^2 \geqslant
1$.

Now recalling the earlier result
\begin{equation}
  \left\| \nabla u \right\|_{L^2} + \left\| \nabla b \right\|_{L^2} \leqslant
  M \left( t \right)^{\delta} ,
\end{equation}
where $\delta$ is given by (\ref{eq:choice-delta}). Such $\delta$ satisfies $A
\delta + B \leqslant 2, \frac{2}{1 + a} \delta + \frac{2}{1 + a} \leqslant 2$.
By denoting $A \left( t \right) \assign 1 + \left\| \mathcal{L}u
\right\|_{L^2}$ and using the facts that $M \left( t \right) > 1$, $\log M
\left( t \right) > 1$, we conclude
\begin{equation}
  \frac{\mathd}{\mathd t} \left( \left\| \nabla^k u \right\|_{L^2}^2 + \left\|
  \nabla^k b \right\|_{L^2}^2 \right) \leqslant CA \left( t \right) M \left( t
  \right)^2 \log \left( M \left( t \right) \right) .
\end{equation}
The integration of this equation, together with the energy inequality, gives
\begin{equation}
  M \left( t \right) \leqslant C \left( T_0 \right)  \left[ 1 + \int_{T_0}^t A
  \left( \tau \right) M \left( \tau \right) \log \left( M \left( \tau \right)
  \right) \mathd \tau \right] .
\end{equation}
Standard Gronwall's inequality then gives
\begin{equation}
  M \left( t \right) \leqslant C \left( T_0 \right)^{\exp \left[ C \left( T_0
  \right)  \int_{T_0}^t A \left( \tau \right) \mathd \tau \right]}
\end{equation}
which is uniformly bounded for all $t \in \left( T_0, T \right)$ since
$\int_{T_0}^T A \left( \tau \right) \mathd \tau < \infty$.

Therefore we have shown that $\left\| u \right\|_{H^k}, \left\| b
\right\|_{H^k}$ are uniformly bounded over $\left( T_0, T \right)$, thus
completing the proof.

\begin{remark}
  \label{rmk:alpha-larger}The case $\alpha > 1 + \frac{n}{2}$ can be proved
  along the same line, with each step much easier. More specifically, in this
  case (\ref{eq:energy-est}) immediately gives $\left\| \nabla u
  \right\|_{L^{\infty}} \in L^2 \left( 0, T \right)$, which leads to a priori
  $H^1$ bounds. This allows us to simply take $\delta = 0$ in the subsequent 
  steps. 
\end{remark}

\appendix\section{Proof of Lemma \ref{lem:log-ineq}}

The proof involves some basic facts from Littlewood-Paley theory, which we
recall here.

Let $\mathcal{S}$ be the Schwartz class of rapidly decreasing functions
and $\widehat{f}(\xi)$ denote the Fourier transform of $f(x)$, i.e.
\begin{equation}
  \widehat{f} \left( \xi \right) \assign \frac{1}{\left( 2 \pi \right)^{n / 2}} 
  \int_{\mathbbm{R}^n} e^{- ix \cdot \xi} f \left( x \right) \mathd x.
\end{equation}
Consider $\phi \in \mathcal{S}$ whose frequency is localized:
\begin{equation}
  \text{\tmop{Supp}}  \widehat{\phi} \subset \left\{ \xi \in \mathbbm{R} :
  \frac{1}{2} \leqslant \left| \xi \right| \leqslant 2 \right\}
\end{equation}
with $\widehat{\phi} \left( \xi \right) > 0$ if $\frac{1}{2} < \left| \xi \right|
< 2$. Now define $\phi_j$ through $\widehat{\phi_j} = \widehat{\phi} \left( 2^{-
j} \xi \right)$. We can multiply $\phi$ by a normalization constant such that
the following holds:
\begin{equation}
  \sum_{j \in \mathbbm{Z}} \widehat{\phi}_j \left( \xi \right) = 1, \hspace{2em}
  \forall \xi \in \mathbbm{R}^n \backslash \left\{ 0 \right\} .
\end{equation}
For any $k \in \mathbbm{Z}$ we can define operators $S_k$ and
$\triangle_k$ by
\begin{eqnarray}
  \widehat{S_k f} \left( \xi \right) & : = & \left[ 1 - \sum_{j \geqslant k +
  1} \widehat{\phi}_j \left( \xi \right) \right]  \widehat{f} \left( \xi \right) \\
  \widehat{\triangle_k f} \left( \xi \right) & \assign & \widehat{\phi}_k \left(
  \xi \right) \widehat{f} \left( \xi \right) . 
\end{eqnarray}
The most important properties of the operators $S_k, \triangle_k$ are the
following Bernstein inequalities: For any $1 \leqslant p \leqslant q \leqslant
\infty$, and $\beta, \beta'$ multi-indices with $\beta \geqslant 0$,
\begin{eqnarray}
  \left\| S_k \partial^{\beta} f \right\|_{L^q} & \leqslant & C 2^{\left|
  \beta \right| k} 2^{kn \left( \frac{1}{p} - \frac{1}{q} \right)}  \left\| f
  \right\|_{L^p} ; \\
  \left\| \triangle_k \partial^{\beta'} f \right\|_{L^q} & \leqslant & C
  2^{\left| \beta' \right| k} 2^{kn \left( \frac{1}{p} - \frac{1}{q} \right)} 
  \left\| f \right\|_{L^p} . 
\end{eqnarray}
Now we are ready to prove Lemma \ref{lem:log-ineq}. The proof is standard and we omit some calculation details. 

{\noindent}{\tmstrong{Proof of Lemma \ref{lem:log-ineq}.}}

We have
\begin{eqnarray}
  \left\| \nabla u \right\|_{L^{\infty}} & \leqslant & \left\| S_{- 1} \nabla
  u \right\|_{L^{\infty}} + \sum_{j = 0}^N \left\| \nabla \triangle_j u
  \right\|_{L^{\infty}} + \sum_{j = N + 1}^{\infty} \left\| \nabla \triangle_j
  u \right\|_{L^{\infty}} \\
  & \leqslant & C \left[ \left\| u \right\|_{L^2} + \sum_{j = 0}^N \frac{2^{j \left(
  1 + \frac{n}{2} \right)}}{g \left( s_j \right)}  \left\| \triangle_j u
  \right\|_{L^2} g \left( s_j \right) + \sum_{j = N + 1}^{\infty} 2^{j \left(
  1 + \frac{n}{2} - k \right)} 2^{kj}  \left\| \triangle_j u \right\|_{L^2} \right]. 
\end{eqnarray}
where $s_j\in \left(2^{j-1},2^{j+1} \right)$ is chosen such that
\begin{equation}
  g \left( s_j \right) \geqslant \frac{1}{2} \sup_{2^{j - 1} < s < 2^{j + 1}}
  g \left( s \right) . \label{eq:sj}
\end{equation}
Now we estimate the second term as follows:
\begin{eqnarray}
  \sum_{j = 0}^N \frac{2^{j \left( 1 + \frac{n}{2} \right)}}{g \left( s_j
  \right)}  \left\| \triangle_j u \right\|_{L^2} g \left( s_j \right) &
  \leqslant & \left[ \sum_{j = 0}^N \left( \frac{2^{j \left( 1 + \frac{n}{2}
  \right)}}{g \left( s_j \right)}  \left\| \triangle_j u \right\|_{L^2}
  \right)^2 \right]^{1 / 2}  \left[ \sum_{j = 0}^N g \left( s_j \right)^2
  \right]^{1 / 2} \\
  & = & C \left[ \sum_{j = 0}^N \left\| \frac{2^{j \left( 1 + \frac{n}{2}
  \right)}}{g \left( s_j \right)}  \widehat{\triangle_j u} \right\|_{L^2}^2
  \right]^{1 / 2}  \left[ \sum_{j = 0}^N \log \left( e + s_j \right)
  \right]^{1 / 2} \\
  & \leqslant & C N \left[ \int_{\mathbbm{R}^n} \left| \frac{\left| \xi
  \right|^{1 + n / 2}}{g \left( \left| \xi \right| \right)} \widehat{u} \left( \xi
  \right) \right|^2 \mathd \xi \right]^{1 / 2} \\
  & = & CN \left\| \frac{\Lambda^{1 + \frac{n}{2}}}{g \left( \Lambda
  \right)} u \right\|_{L^2} . 
\end{eqnarray}
Here we have used the assumption (\ref{eq:cond-main}), the definition of $s_j$
(\ref{eq:sj}), the Plancherel theorem, and the following facts about
$\widehat{\phi}_j \left( \xi \right)$: 1. $\tmop{supp} \left( \widehat{\phi}_j \right)
\subseteq \left\{ 2^{j - 1} < \left| \xi \right| < 2^{j + 1} \right\}$; 2. $0
\leqslant \widehat{\phi}_j \left( \xi \right) \leqslant 1 \Longrightarrow \left|
\widehat{\phi}_j \left( \xi \right) \right|^2 \leqslant \widehat{\phi}_j \left( \xi
\right)$; 3. $\sum_{j = 0}^N \widehat{\phi}_j \left( \xi \right) \leqslant 1$.

For the third term we have
\begin{eqnarray}
  \sum_{j = N + 1}^{\infty} 2^{j \left( 1 + \frac{n}{2} - k \right)} 2^{kj} 
  \left\| \triangle_j u \right\|_{L^2} 
   & \leqslant & \left[ \sum_{j = N +
  1}^{\infty} 2^{2 j \left( 1 + \frac{n}{2} - k \right)} \right]^{1 / 2}   \left[ \sum_{j = N + 1}^{\infty} 2^{2 kj}  \left\| \triangle_j u
  \right\|_{L^2}^2 \right]^{1 / 2} \\
  & \leqslant & 
  2^{\left( 1 + \frac{n}{2} - k \right) N}  \left\| u
  \right\|_{H^k} . 
\end{eqnarray}
Summarizing, we have
\begin{equation}
  \left\| \nabla u \right\|_{L^{\infty}} \leqslant C \left[ \left\| u \right\|_{L^2}
  + N \left\| \frac{\Lambda^{1 + \frac{n}{2}}}{g \left( \Lambda \right)} u
  \right\|_{L^2} + 2^{\left( 1 + \frac{n}{2} - k \right) N}  \left\| u
  \right\|_{H^k}\right] .
\end{equation}
Taking $N$ such that $2^{\left( k - 1 - \frac{n}{2} \right) N} \approx \left\|
u \right\|_{H^k}$ gives the result.

\paragraph{\textbf{\\ Acknowledgment}}X. Yu and Z. Zhai are supported by a grant from
NSERC and the Startup grant from Faculty of Science of University of Alberta.
The authors would like to thank the anonymous referees for their valuable comments and suggestions.

\end{document}